\documentclass{article}

\usepackage[utf8]{inputenc}
\usepackage{authblk}
\usepackage{amssymb}
\usepackage[pdftex]{graphicx}
\usepackage{epstopdf}
\graphicspath{ {./images/} } 
\usepackage{amsmath,textcomp}
\usepackage{bbm}
\def\1{\mathbbm{1}}

\DeclareMathOperator\erf{erf}
\DeclareMathOperator{\spann}{span}

\title{\bf On the basic reproduction number in continuously structured populations}
\author[*]{Carles Barril}
\author[*]{\`Angel Calsina}
\author[*]{S\'{i}lvia Cuadrado}
\author[**]{Jordi Ripoll}
\affil[*]{\it Departament de Matem\`atiques, Universitat Aut\`onoma de Barcelona}
\affil[**]{\it Departament d'Inform\`atica, Matem\`atica Aplicada i Estad\'\i stica, Universitat de Girona}
\affil[ ]{\small \texttt{carlesbarril@mat.uab.cat, acalsina@mat.uab.cat, silvia@mat.uab.cat, jripoll@imae.udg.edu}}

\begin{document}

\maketitle


\begin{abstract}
%
%
In the framework of population dynamics, the basic reproduction number $\mathcal{R}_0$ is, by definition, the expected number of offspring that an individual has during its lifetime. In constant and time periodic environments it is calculated as the spectral radius of the so-called \textit{next-generation operator} (\cite{DHM,inaba2}). In continuously structured populations defined in a Banach lattice $X$ with concentrated states at birth one cannot define the next-generation operator in $X$. In the present paper we present an approach to compute the basic reproduction number of such models as the limit of the basic reproduction number of a sequence of models for which $\mathcal{R}_0$ can be computed as the spectral radius of the next-generation operator. We apply these results to some examples: the (classical) size-dependent model, a size structured cell population model, a size structured model with diffusion in structure space (under some particular assumptions) and a (physiological) age-structured model with diffusion in structure space.

\end{abstract}
\textbf{Key words}: next generation operator; basic reproduction number ; physiologically structured population.
%

%
%

\section{Introduction}
In a deterministic model of population growth, the basic reproduction number $\mathcal{R}_0$ is defined as the expected number of offspring that an individual has throughout his life (in an epidemiological model, as the expected number of new infections a newly infected individual will produce). The main interest of the basic reproduction number is that, as it is intuitively clear, a small population will (begin to) establish in a given environment if $\mathcal{R}_0  > 1$ whereas it will become extinct whenever $\mathcal{R}_0  < 1$. Similarly, provided some infected individual is present, an epidemic outbreak will arise if $\mathcal{R}_0  > 1$ whereas the opposite strict inequality guarantees that the epidemic will not occur.

In structured populations, when the birth event can happen in different individual states (of size, phenotype, spatial position, etc.) one talks of “typical” individual but it is not always clear what “typical” individual means and so, which is the expected number of offspring of a ”typical” individual.

On the other hand, in deterministic modeling one could be concerned with a definition that involves a probability concept as expectation and would like to have a deterministic definition as well.

Both drawbacks are solved by defining $\mathcal{R}_0$ as the spectral radius of the so-called \textit{next-generation operator}, which maps a (distribution of) population to the (distribution of) population of their offspring along the whole life span of the former (\cite{DHR}, \cite{Thieme}, \cite{DHM}, \cite{inaba2}).

Nevertheless, in some models the definition of the birth operator and that of the next-generator operator depend on the choice of the sometimes arbitrary concept of birth event (\cite{CD}, \cite{ripoll1}).
Hence more than one basic reproduction number can be meaningfully defined in some models (\cite{CD}, \cite{ripoll1}).

When the population dynamics can be described by a linear system of odes, the next-generation operator is a matrix for which the spectral radius can be computed and therefore there is no problem with the definition above leaving apart that it may be non-uniquely determined, as commented in the previous paragraph. Fortunately (we should better say inevitably) the sign of  $\mathcal{R}_0-1$ coincides with that of the real part of the so-called Malthus constant (the first eigenvalue of the matrix defining the system of odes), independently of the choice of the decomposition of the latter in a mortality/transition operator and a birth operator provided some natural assumptions hold (\cite{DHR}, \cite{CD}). The sign relation also holds in infinite dimension \cite{DHM}.

However, in continuously structured populations with concentrated state at birth, the birth rate shows as a boundary condition instead of a birth operator. This makes it difficult (or impossible) to obtain the next-generation operator and thus adapt to these type of models in  infinite dimensional spaces (which we will call type II models) the above definition of $\mathcal{R}_0$. In order to overcome this difficulty, in the present paper we consider some examples where $\mathcal{R}_0$ is defined for type II models as a limit of basic reproduction numbers of approximate models with distributed states at births for which the next generation operator can be defined as a bounded linear operator, which we call type I models.

More explicitly, we will consider diffusion-convection partial differential equations of the form of conservation laws for the density of individuals with respect to some one-dimensional continuous variable which stands for a physiological state as size or age or for an external state as spatial position, with non-flux boundary conditions and with an incoming flow of new individuals via a birth operator whose range is the span of a function in the space of states giving the distribution of the offspring (for which the next generation operator is well defined).

Next we will assume a sequence of type I models such that this distribution of offspring concentrates at some point in the closure of the domain of the structuring variable. This sequence tends to a conservation law where the birth term appears as a boundary condition and consequently the next generation operator is not definable as a bounded linear operator in the state space. Many of the classic age and size structured population models are of this second class which we call type II models. For this we define its basic reproduction number as the limit (it exists uniquely under suitable hypotheses) of the basic reproduction numbers of the sequence of type I models tending to it.

A different approach to overcome the difficulty to define the next generation operator in the case of concentrated states at birth is to consider the space of measures as state space. This is done for discrete-time population dynamics in the recent reference \cite{Thieme2}.

The paper is organized as follows: in Section 2 we give, in a general framework, the definition of $\mathcal{R}_{0}$ for type II models as the limit of the sequence of basic reproduction numbers of type I models given by conservation laws for the density of individuals with respect to some one-dimensional continuous variable. Section $3$ is devoted to the application of the results in Section 2 to some examples: the (classical) size-dependent model, a size structured cell population model, a size structured model with diffusion in structure space (under some particular assumptions) and a (physiological) age-structured model with diffusion in structure space. For this last example, the computation of the basic reproduction number for the sequence of type I models is done in two different ways, one of them included in the Appendix. At the end of Section 3, a brief discussion on the dependence of $\mathcal{R}_{0}$ on the diffusion coefficient is undertaken with the conclusion that, for some age specific birth functions, there is an optimal value of the diffusion coefficient which maximizes $\mathcal{R}_{0}$. Therefore, the presence of a moderate diffusion in age tends to increase $\mathcal{R}_0$ and hence the survival chances of a population.

\section{General framework}\label{general}
 Some models in population dynamics (which we call type I models, see \cite{ripoll1}) can be described by non-linear abstract differential equations for the density of individuals with respect to some structuring variables. Specifically, let $X$ be a Banach lattice and let $u(t) \in X$ be the distribution of individuals at time $t\geq 0$, then
the original non-linear ecological problem can be decomposed as
\begin{equation} \label{Non-linear}
u'(t)= B(u(t)) u(t) - M(u(t)) u(t)\, , \, u(0)=u_0 \in X
\end{equation}
where for any $u$, $B(u)$ is the linear birth operator and for any $u$, $M(u)$ is the linear operator corresponding to non-birth terms (mortality and transitions in general). Of course, as mentioned before, this decomposition depends on what is considered as a birth event (\cite{CD, ripoll2}). For the study of possible extinction or settlement of the population we can focus on the so-called extinction steady-state $u^*=0$ and its stability which is analyzed by means of the following formal linearized model
$$
u'(t)= B u(t) - M u(t)\, , \, u(0)=u_0 \in X \, .
$$
where we set $B(0):=B$ and $M(0):=M$ with some abuse of notation.

In order to proceed, we need the following assumptions: $B: \mathcal{D}_{B} \subset X \to X$ is a positive linear operator, and $M: \mathcal{D}_M \subset X \to X$ is such that $-M$ is the generator of a strongly continuous semigroup of positive linear operators $T(t):= e^{-M t}$, whose spectral bound is strictly negative $s(-M)<0$ (so that the population goes to extinction when there are no births). Therefore, 0 belongs to the resolvent set of $M$ and there exists $\int_0^\infty T(t)\, dt=M^{-1} $. Moreover we assume that $\mathcal{D}_{M} \subset \mathcal{D}_{B}$ and that the linear operator  $BM^{-1}$ is bounded. That is guaranteed for instance when $B$ is bounded but this is not a necessary condition (see \cite{ripoll2}, Sect. 5). Finally we assume that $B-M$ is the infinitesimal generator of a positive semigroup.

The operator $BM^{-1}$ can be interpreted as the next generation operator. For the finite dimensional case, see for instance \cite{CD}. In the infinite dimensional case, if $B$ is bounded see \cite{DHM} and Section 3 in \cite{ripoll2}. If $B$ is not bounded, $BM^{-1}$ is still interpreted as the next-generation operator, see \cite{Thieme} and Section 5 in \cite{ripoll2}.

For type I models, we have that the next-generation operator $BM^{-1}$ is defined on the same Banach space $X$ as the population density. The basic reproduction number $\mathcal{R}_0$ (which is an alternative to the Malthusian parameter) is then given by the spectral radius of the positive bounded linear operator $BM^{-1}$ (\cite{DHM,ripoll2}). Moreover, $\mathcal{R}_0$ is a non-negative spectral value which is actually the largest $\lambda \geq 0$ for which the operator $B - \lambda M$ does not have a bounded inverse, see e.g. \cite{ripoll1}.

Often $\mathcal{R}_0$ cannot be computed explicitly and then the following upper bound is useful to find a sufficient condition for extinction
 $$\displaystyle\mathcal{R}_0\leq ||BM^{-1}||_X= \sup\left\{ \frac{||BM^{-1}\phi||_X}{||\phi||_X} : \phi \in X\setminus\{0\} \right\} =$$
$$
=\sup\left\{ \frac{||B\psi||_X}{||M\psi||_X} : \psi \in \mathcal{D}_M\setminus\{0\} \right\}.
$$

It is important to notice that, although type I models contain a large family of population models (in particular, those of discrete structure, for which $\dim X<\infty$), often some structured models cannot be cast into this form.

Indeed, beyond the "standard" case considered above where $\mathcal{R}_0 = \rho ( B M^{-1})$, we can consider another
family of population models (which we call type II models as in \cite{ripoll1}) where the
population distribution at time $t\geq 0$ is such that $u(t) \in X$
and the distribution of newborns per
time-unit does not belong to $X$,
and so the model cannot be written in the form (\ref{Non-linear}). For instance, in the
case of age-structured populations, all newborns are concentrated at age $0$ which translates as a boundary
condition for an equation of the form $u'(t)= - M(u(t)) u(t)$. More in general this is very often the
case when the individual state space, i.e., the set where the
structuring variable lives (its domain) is a continuum, for instance
a real interval, and the set of states at birth is discrete.

In the forthcoming we will consider models given by
conservation laws for densities of individuals living on an interval
$I$, where often the newborns show as a flow crossing the boundary
of this domain.

Notice that in the case of type II models one
cannot, in principle, define a next generation operator.

In order to define $\mathcal{R}_{0}$ for type II models, following the approach developed in \cite{ripoll1}, we will consider a specially-built sequence of type I
models in $X = L^1(I)$
\begin{equation}\label{type1}
u'(t)= B_k u(t) - M u(t)\, , \, u(0)=u_0 \in X, \, \, k\in \mathbb{N}
\, .
\end{equation}
The linear operator $M: \mathcal{D}_M \subset X \to X$,
as before, is such that $-M$ is the generator of a
strongly continuous semigroup of positive linear operators whose
spectral bound is strictly negative $s(-M)<0$ and $B_k: \mathcal{D}_{B_{k}}\subset X \to X$  are positive linear operators, all of
finite rank $m<\infty$ such that $B_{k}M^{-1}$ is bounded and $B_{k}-M$ is the generator of a strongly continuous semigroup. To fix ideas let us
assume that $m=1.$ We further assume that, for any $k$, $B_k u = (Lu)\,
\varphi_k$ for some positive linear functional $L$  and some positive $\varphi_k \in X.$ Notice that $\varphi_k$ represents the distribution of offspring.
Moreover we assume that the sequence $M^{-1} \varphi_k$ converges in $X$
to a certain function $\psi_{\infty} \in \mathcal{D}_{L}=\mathcal{D}_{B_{k}}$, and that $LM^{-1}\varphi_k$ converges to $L \psi_{\infty}$. Notice that the last hypothesis (and the fact that $\psi_{\infty} \in \mathcal{D}_{L}$) authomatically hold if $B$ is bounded.

Then we will have that the
sequence of basic reproduction numbers corresponding to the
approximate models has a limit, which we claim can be named as the
basic reproduction number of the limit type II model:
\begin{equation}\label{R0}
\mathcal{R}_0 := \lim_{k \to \infty} \mathcal{R}_{0,k} = \lim_{k \to
\infty} \rho ( B_k M^{-1})= \lim_{k \to \infty} \rho ((L M^{-1}
\cdot)\varphi_k) = \lim_{k \to \infty} L M^{-1} \varphi_k = L
\psi_{\infty} \, ,
\end{equation}
where $\rho(\cdot )$ stands for the spectral radius of the positive
bounded linear operators involved and where we have used that $range
(B_k) = \spann \{\varphi_k\}$.\\
Of course this definition needs a proof of correctness in the sense
that it does not depend on the sequence ${\varphi_k}$ chosen to
approximate the limit model. More precisely, we have to clarify what
we do understand by limit model of a sequence of type I models of
the form (\ref{type1}) with $B_k u = (Lu)\, \varphi_k$ and check
that two sequences of functions ${\varphi_k}$ and
${\hat{\varphi}_k}$ giving rise to the same limit model are such
that ${M^{-1}\varphi_k}$ and
${M^{-1}\hat{\varphi}_k}$ tend to the same limit.\\
Thus let us consider a sequence of conservation laws in the form
\begin{equation}\label{conservation}
\frac{d}{dt} \int_{a}^{b} u(s,t) ds = \Phi(a,t) - \Phi(b,t)
-\int_{a}^{b} \mu(s)u(s,t) dt + Lu(\cdot, t) \int_a^b
\varphi_k(s) ds
\end{equation}
for the rate of change of the population number of individuals whose individual state (size, spatial position, age, etc.) belongs to the
 (arbitrary) interval $(a,b) \subset I.$ Here $\Phi(x,t)$ is the flow crossing the point $x$ at time $t$ and we assume that it is given by a linear combination
 of the density $u$ and its partial derivative $u_s,$ $\mu$ is the state-dependent per capita death rate (so the integral above gives the rate of loss of individuals
 with individual state within the interval $[a,b]$). We also assume that the new individuals are being born according to a fixed distribution of
 individual states $\varphi_k \in X,$ with $\int_{I}\varphi_k (s) ds = 1.$ Hence the last term stands for the population birth rate.\\
This leads, for a smooth density $u,$ to the following sequence of partial
differential equations
\begin{equation}\label{type1pde}
\partial_t u(x,t) = - \partial_x \big(\gamma(x)u(x,t) - D(x) \partial_x u(x,t)\big) - \mu(x) u(x,t) + Lu(\cdot, t)\, \varphi_k(x)
\end{equation}
(where $D(x) > 0$ or $D(x) \equiv 0)$, with non-flux boundary conditions at the (possibly empty) boundary
of $I.$ In \eqref{type1pde} we can firstly identify an unbounded linear operator, formally
written
$$(Mu)(x) =
\partial_x \big(\gamma(x)u(x) - D(x) \partial_x u(x)\big) + \mu(x) u(x),
$$
with non-flux boundary conditions, for transition between states and mortality and secondly rank one linear operators $(B_{k}u) (x)= Lu \, \varphi_k(x)$
corresponding to reproduction and verify that the model
is of the type defined in (\ref{type1}).\\
Let us further assume that the sequence of states distributions at
birth ${\varphi_k}$ tends to concentrate, for instance around a
point in the closure of $I$. More precisely, that there exists an
individual state $x_0 \in \overline{I}$ such that for any $[a,b] \subset \overline{I}$, $\lim_{k \rightarrow
\infty} \int_a^b \varphi_k (s) ds = 1$ whenever $x_0 \in [a,b]$ and
$\lim_{k \rightarrow \infty} \int_a^b \varphi_k (s) ds = 0$
otherwise.\\
Assuming the second case (i.e. $x_0 \notin [a,b]$) and making $k$ go to infinity in
(\ref{conservation}) we simply obtain the partial differential
equation
\begin{equation}\label{pde}
\partial_t u(x,t) = -
\partial_x \big(\gamma(x)u(x,t) - D(x) \partial_x u(x,t)\big) - \mu(x) u(x,t)
\end{equation}
for $x \neq x_0$.

On the other hand, if $[a,b]$ contains $x_0,$ making $k$ go to
infinity in (\ref{conservation}) we first obtain the integral
equation
$$
\frac{d}{dt} \int_{a}^{b} u(s,t) ds = \Phi(a,t) - \Phi(b,t)
-\int_{a}^{b} \mu(s)u(s,t) dt + Lu(\cdot, t)
$$
and afterwards, making that the interval $[a,b]$ collapses to $x_0,$ we obtain the
condition
\begin{equation}\label{bcflux}
0 = \lim_{a \rightarrow x_0^{-}}\Phi(a,t) - \lim_{b \rightarrow
x_0^{+}}\Phi(b,t) + Lu(\cdot, t).
\end{equation}
This can be seen as a boundary condition when $x_0$ belongs to the
boundary of $I.$ For instance if $x_0$ is the left endpoint of $I$
we necessarily have $a = x_0$ and $\lim_{a \rightarrow
x_0^{-}}\Phi(a,t) = \Phi(x_0,t) = 0$ for the non-flux boundary condition and (\ref{bcflux}) reduces to a
(in general) Robin type boundary condition
\begin{equation}\label{bc1}
 \lim_{b \rightarrow x_0^{+}}\Phi(b,t) = \gamma(x_0)u(x_0,t) - D(x_0) \partial_x u(x_0,t) = Lu(\cdot, t).
\end{equation}

Notice, as particular cases, that the latter supplementing the pde
(\ref{pde}) with $D(x)\equiv 0$  corresponds to the well known
equation of the linear size dependent population dynamics when individuals
are all born with the same size (and in particular to the equation of the linear age
dependent population dynamics if, furthermore, $\gamma(x)\equiv
1)$ whereas when $D(x)$ is strictly positive the pde \eqref{pde} could be interpreted as a
linear convection-diffusion equation (purely diffusive if $\gamma(x)
\equiv 0)$ for a population living in a one-dimensional habitat such
that, in order to reproduce, adult individuals travel instantly to
the point $x_0$ (a common ``nest") and instantly go back after, to
where they were living.\\

On the other hand, if $x_0$ belongs to the interior of $I$,
(\ref{bcflux}) causes a discontinuity in the flux at the point
$x_0.$ Biologically this could correspond to a situation like the
latter, but where the "nest" is not at the boundary of the domain.
From the point of view of the pde (\ref{pde}) one can consider it on
two subdomains of $I$ separated by $x_0$ and with "boundary
conditions" at $x_0$ given by
\begin{equation}\label{bc2}
\lim_{x \rightarrow x_0^{+}} \big( \gamma(x)u(x,t) - D(x) \partial_x
u(x,t)\big) - \lim_{x \rightarrow x_0^{-}} \big(\gamma(x)u(x,t) -
D(x)
\partial_x u(x,t) \big) = Lu(\cdot, t).
\end{equation}
As one could expect, we define the limit model of a sequence of type
I models of the form (\ref{type1pde}) (a particular case of the
abstract form (\ref{type1})) as the one given by the partial
differential equation (\ref{pde}) plus boundary conditions of the
form (\ref{bc1}) or (\ref{bc2}).\\
Next we show that different sequences $\varphi_k$ leading to the same
limit model (i.e. concentrating at the same point $x_0$) share the
same limit of the sequence $M^{-1}\varphi_k$ and thus that the basic
reproduction number of the limit model (given by (\ref{R0})) is well
defined. Indeed, let us denote by $G(x,s)$ the Green's function of the
operator $M$, i.e., such that one can write
$$
(M^{-1}\varphi)(x) = \int_{I} G(x,s) \varphi(s) ds,
$$
(for a comprehensive treatise on Green's functions one can see \cite{Stakgold}).

Then, due to the hypotheses on $\varphi_k,$
$$
\psi_{\infty}(x) = \lim_{k\rightarrow \infty}(M^{-1}\varphi_k)(x) =
\lim_{k\rightarrow \infty} \int_{I} G(x,s) \varphi_k(s) ds =
G(x,x_0),
$$
independently of the choice of the sequence. Moreover we
can now write
\begin{equation}
\label{green}
\mathcal{R}_0 := \lim_{k \to \infty} \mathcal{R}_{0,k} = \lim_{k \to \infty}LM^{-1}\varphi_k= L
\psi_{\infty} = L G(\cdot,x_0).
\end{equation}

\section{Examples}
\subsection{A size dependent model}
\label{size}
As a first example let us consider the type II model associated to the classical (linear) size dependent population
dynamics with potentially unbounded size. Then, (\ref{pde}) reduces
to
$$
\partial_t u(x,t) = -
\partial_x \big(\gamma(x)u(x,t)\big) - \mu(x) u(x,t)=: (Mu)(x,t) \quad x \in (x_0, +\infty)
$$
for $u(\cdot,t) \in L^{1}(x_0, \infty)$, where $x_0$ denotes the individual's size at birth and with boundary condition
$$\gamma(x_0)u(x_0,t) = \int_{x_0}^{\infty} \beta (s) u(s,t)ds=:(Lu)(t)$$
 (a particular form of (\ref{bc1})). We assume that the individual growth rate $\gamma(x)$ and the mortality rate $\mu(x)$ are smooth and bounded below by a positive number and that $\beta(x)$ is bounded and such that $\beta(x) \geq 0$. Moreover, to prevent that the individuals reach infinite size in finite time we assume $\int_{x_0}^{\infty}\frac{1}{\gamma(x)}\mbox{d}x=\infty.$\\
We can choose a sequence $\varphi_k(x)$ concentrating at $x_0$ and consider this model as
a limit of type I models of the form (\ref{type1pde}) with $D(x) \equiv 0$. Notice that $-M$ and $B_{k}-M$ (where $B_{k}u=(Lu)\varphi_{k}$) are both generators of strongly continuous semigroups. \\ Using the variation of constants formula it is straightforward
to compute
$$
\begin{array}{rcl}
\big(M^{-1}v\big)(x) &= &\int_{x_0}^{x} v(s)e^{\int_{x_0}^s
\frac{\mu(\sigma)}{\gamma(\sigma)} d\sigma} ds \,\frac{1}{\gamma(x)}e^{-\int_{x_0}^x \frac{\mu(s)}{\gamma(s)} ds}\\ \\ & =& \int_{x_0}^{\infty}\frac{1}{\gamma(x)}e^{-\int_s^x \frac{\mu(\sigma)}{\gamma(\sigma)} d\sigma}H(x-s)v(s)ds,
\end{array}
$$
where $H$ is the Heaviside function. Hence, $G(x,s) = \frac{1}{\gamma(x)}e^{-\int_s^x \frac{\mu(\sigma)}{\gamma(\sigma)} d\sigma}\, H(x-s),$ and we have
$$
\psi_{\infty} (x) = G(x,x_0) = \frac{1}{\gamma(x)}e^{-\int_{x_0}^x \frac{\mu(\sigma)}{\gamma(\sigma)} d\sigma}.
$$
Using (\ref{green}) we obtain
\begin{equation}
\label{brnsize}
\mathcal{R}_0 = L \psi_{\infty} = \int_{x_0}^{\infty}
\frac{\beta(x)}{\gamma(x)}e^{-\int_{x_0}^x \frac{\mu(s)}{\gamma(s)} ds} dx
\end{equation}
which coincides with the basic reproduction number in size
dependent population dynamics obtained in the literature (\cite{IP}, \cite{inaba2}) as the number of newborns that an individual is expected to produce during his reproductive life. \\

\subsection{A size-structured cell population model}

Let us start by considering a cell population structured by size $x$ and
such that cells divide when they reach a given size (normalized to
$1$) giving rise to two daughter cells whose size is not necessarily
the same. Indeed, let us assume $x \in (0,1)$ and a size
distribution of the newborn cells given by a probability density function $\varphi(x).$
Conservation of the mass during cell division forces that
$\varphi(x) = \varphi (1-x),$ i.e. the distribution of newborn cells
is symmetric with respect to half the mother cell size. In this way
we obtain a type I model given by the pde
\begin{equation}
\label{cell}
\partial_t u(x,t) = -\partial_x(\gamma(x) u(x,t)) - \mu(x) u(x,t) + 2 \gamma(1) u(1,t) \varphi(x), \,\, x \in \big(0,1 \big)
\end{equation}
with zero flux at $x = 0$ for the density of cells $u(x,t)$ and
where $\gamma$ is the individual growth rate, $\mu$ the mortality rate
and the last term is the (distributed) birth rate. We assume the same hypotheses on $\gamma$ and $\mu$ of the previous section. The pde is of the
form \eqref{type1pde} with, as before $D(x) \equiv 0$.

Notice that in this case the birth operator $(B u) (x) = (L u)
\varphi(x) = 2 \gamma(1) u(1) \varphi(x)$, despite being of rank $1$, is
only relatively bounded with respect to the linear operator $(-M u)
(x) = -\partial_x(\gamma(x) u(x)) - \mu(x) u(x)$ with non-flux boundary
condition at the left hand end $x_0 = 0$ (i.e. $B$  is not bounded as a linear operator in $L^1$ but $B M^{-1}$ is bounded, which amounts to that $LM^{-1}$ is a bounded linear form). \\
In order to be in the hypotheses of Section \ref{general} we need $B-M$ to be the infinitesimal generator of a strongly continous positive semigroup. The fact that $B$ is relatively bounded is in general not enough to ensure that $B-M$ is a generator of a strongly continuous semigroup \cite{Engel}. Nevertheless it suffices that $B$ is a bounded linear operator in the domain of $(-M)$ endowed with the graph norm (see \cite{Engel}, Cor. III.1.5).\\
For \eqref{cell} this is satisfied provided that  $\varphi \in D(-M) = \{u \in W^{1,1}(0,1): u(0)=0\}.$  Indeed, we have for $u \in  \mathcal{D}_{M},$ and using that $0$ belongs to the resolvent set of the operator $M$,
$$
||Bu||_{\mathcal{D}_{M}}=||MBu|| =
||M\, (Lu) \varphi|| = \mid Lu \mid ||M \varphi ||= \mid LM^{-1}Mu \mid ||M \varphi ||
$$
$$
\leq  ||M \varphi ||\,||L^{-1}M||\, ||M u|| = C
||u||_{\mathcal{D}_{M}},
$$
where the norms without subscript are $L^1-$norms.
\\

Since the mortality/transition operator $M$ is exactly the same as in the previous section, the Green's function is
$$
G(x,s) = \frac{e^{-\int_s^x \frac{\mu(\sigma)}{\gamma(\sigma) } d
\sigma}}{\gamma(x)} H(x-s).
$$

So, the basic reproduction number for this model is given by
$$
\mathcal{R}_0 = L M^{-1} \varphi =2 \gamma(1)(M^{-1}\varphi)(1) = 2 \gamma(1) \int_0^1 G(1,s) \varphi(s)
ds = 2 \int_0^1 e^{- \int_s^1 \frac{\mu(\sigma)}{\gamma(\sigma)} d
\sigma} \varphi (s) ds.
$$

A bit less realistic but apparently simpler, one can
alternatively assume that the division is perfectly symmetric. To
deal with this assumption, similarly to what we have done in the previous example, we can
choose a sequence of functions $\varphi_k \in D(-M)$ concentrating at $x_0 = 1/2$ (for
instance $\varphi_k (x) = a_k x (1-x) e^{-k \mid
x- \frac{1}{2}\mid}$ with $a_k$ such that $\int_0^1 \varphi_k = 1$) and consider the model with symmetric division as a limit of a
sequence of type I models as the ones described above (of the
form \eqref{type1pde}).\\
As limiting equation we get
\begin{equation}\label{pdeex}
\partial_t u(x,t) = -\partial_x(\gamma(x) u(x,t)) - \mu(x) u(x,t), \,\,
\end{equation}
with "boundary condition" (according to \eqref{bc2})
\begin{equation}\label{bcex}
\gamma(\frac{1}{2}) \bigg(\lim_{x \rightarrow (\frac{1}{2})^{+}}u(x,t)-
\lim_{x \rightarrow (\frac{1}{2})^{-}}u(x,t)\bigg) = L u(\cdot,t) =
2 \gamma(1) u(1,t),
\end{equation}
which can be considered a birth term giving the influx of new cells
(double of the number of cells that disappear because they
divide).\\

We have that
$$
M^{-1} \varphi_k = \int_0^{1} G( \cdot,s) \varphi_k(s) ds
$$
tends to $G(\cdot,\frac{1}{2})$ in $L^1(0,1)$ and uniformly on $(0,1) \smallsetminus (\frac{1}{2}-\varepsilon, \frac{1}{2}+ \varepsilon)$ (for any $\varepsilon < \frac{1}{2}$).

Therefore, the basic reproduction number for the limit model can be
computed as
$$
\mathcal{R}_0 = \lim \mathcal{R}_{0,k} = \lim LM^{-1} \varphi_k =\lim 2 \gamma(1) \int_0^1 G(1,s) \varphi_k(s)
ds = 2 \gamma(1) G\big(1,\frac{1}{2}\big) = 2 e^{-\int_{\frac{1}{2}}^1 \frac{\mu(\sigma)}{\gamma(\sigma) } d \sigma}.
$$

Finally notice that, for any initial condition, the solution of the
initial value problem for \eqref{pdeex}-\eqref{bcex} will vanish
identically on the interval $(0,\frac{1}{2})$ after a finite time
which means that we can consider equivalently the pde \eqref{pdeex}
on the interval $(\frac{1}/{2},1)$ with boundary condition $\gamma(1/2) u(1/2,t) =
2 \gamma(1) u(1,t).$

\subsection{Size-structured particular model with diffusion}

Let us consider a closed population living in a specific habitat where individuals are classified according to some biometric measure (physiological size) which we assume to evolve in a diffusive way, i.e. random fluctuation in the individual growth. For some examples on structured population models in structure space see for instance \cite{Hadeler}, \cite{PM}, \cite{FH}.

Let $X= L^1(x_0, \infty)$, $x_0\geq 0$, and let $u(\cdot, t)\in X$ be the
density with respect to size at time $t\geq 0$. The non-linear problem is described as:
$$\left\{\begin{array}{l}
 \partial_t u(x,t) + \partial_x \big[\gamma(x,S(t)) u(x,t)-D \partial_x u(x,t)\big]+ \mu(x,S(t)) u(x,t) = 0 \medskip\\
 \gamma(x_0,S(t)) u(x_0,t) - D \partial_x u(x_0,t) = \int_{x_0}^\infty \beta(x,S(t)) u(x,t)\, dx \medskip \\
 S(t)= \int_{x_0}^\infty \sigma(x) u(x,t)\, dx
\end{array}\right.
$$
where $S(t)$ is a weighted population size with $\sigma(x) \geq 0$, $\gamma(x,S(t)) > 0$ is the individual growth rate, i.e. $\frac{dx}{dt}=\gamma(x,S(t))$, $x(0)=x_0$, $D > 0$ is the diffusion coefficient, $\beta(x,S(t))\geq 0$ is the fertility rate and $\mu(x,S(t))\geq 0$ is the mortality rate.
Notice that all newborns are assumed to have size $x_0\geq 0$, and \textit{growth and vital} processes here are density-dependent accounting for limited resources (crowding effects) and size-specific (heterogeneity of the population).
Finally, notice that we would get the classical size-structured model if we had $D \equiv 0$.

The linearization around the trivial steady-state is given by:
\begin{equation}
\label{linearized}\left\{\begin{array}{l}
 \partial_t u(x,t) + \partial_x \big[\gamma(x) u(x,t)-D \partial_x u(x,t)\big]+ \mu(x) u(x,t) = 0 \medskip\\
 \gamma(x_0) u(x_0,t) - D \partial_x u(x_0,t) = \int_{x_0}^\infty \beta(x) u(x,t)\, dx
\end{array}\right.
\end{equation}
where we set $\gamma(x):=\gamma(x,0)$, $\beta(x):=\beta(x,0)$ and $\mu(x):=\mu(x,0)$ with an abuse of notation. We assume the same hypotheses on $\gamma(x)$, $\mu(x)$ and $\beta(x)$ as in Section \ref{size} .

We can choose $\varphi_{k}(x)= k  \mathbbm{1}_{[x_0,x_0+1/k]}(x)$ and consider model \eqref{linearized} as a limit of
the following sequence of type I models in $X= L^1(x_0, \infty)$:

$$\left\{\begin{array}{l}
 \partial_t u(x,t) + \partial_x \big[\gamma(x) u(x,t)-D \partial_x u(x,t)\big]+ \mu(x) u(x,t) =
 \int_{x_0}^\infty \beta(x) u(x,t) dx \cdot k  \mathbbm{1}_{[x_0,x_0+1/k]}(x) \medskip\\
 \gamma(x_0) u(x_0,t) - D\partial_x u(x_0,t) =  0\; ,  \qquad k\in \mathbb{N} \; .
\end{array}\right.
$$
Here, the birth operators are rank one bounded linear operators of
the form $B_k: X \to X$,
$$\big(B_k \phi\big)(x)= \int_{x_0}^\infty \beta(x) \phi(x) dx \cdot
k\, \mathbbm{1}_{[x_0,x_0+1/k]}(x)$$ and the transition/mortality
operator is given by the unbounded linear operator $M: \mathcal{D}_M
\subset X \to X$,
$$\big(M\phi\big)(x)=  (\gamma(x) \phi(x) -D \phi'(x))' + \mu(x) \phi(x)$$ with domain
$\mathcal{D}_M=\left\{\phi \in X : \phi', (\gamma \phi - D \phi')' \in X, \textrm{and } \gamma(x_0) \phi(x_0)-D \phi'(x_0)=0 \right\}$. 

The fact that $-M$ is the generator of a positive semigroup goes back to the work by Feller, Hille and Yosida (see \cite{feller} and the references therein). Since $B_{k}$ is bounded, $B_{k}-M$ is the infinitesimal generator of a positive semigroup.

 For each $k\geq 1$, $\mathcal{R}_{0,k}$ is computed as the largest $\lambda \geq 0$ for which the operator $B_{k}-\lambda M$ does not have a bounded inverse, that is, we have to study the problem $B_k\psi - \lambda M \psi= \xi$, for $\psi\in \mathcal{D}, \xi\in X$, which is,
$$
 \begin{array}{c}
\int_{x_0}^\infty \beta(x) \psi(x) dx \cdot k \, \mathbbm{1}_{[x_0,x_0+1/k]}(x)-  \lambda \big[(\gamma(x) \psi(x) -D \psi'(x))' + \mu(x) \psi(x)\big]= \xi(x) \medskip \\
\gamma(x_0) \psi(x_0)-D \psi'(x_0)=0\end{array}
$$
and integrating the first equation along the whole size-span we get
$$
\int_{x_0}^\infty \big[ \beta(x) -  \lambda \mu(x)\big]\psi(x)\, dx =\int_{x_0}^\infty \xi(x)\, dx .
$$
There is a special case in which the basic reproduction number has an explicit expression. Namely,
the ``\textit{proportional vital rates case}'':
 \begin{equation}\label{birth}\beta(x)= \tilde{\beta} \cdot \mu(x)\, , \; \tilde{\beta}\geq 0 \end{equation} where the relation above becomes $\big[ \tilde{\beta} -  \lambda \big] \int_{x_0}^\infty \mu(x) \psi(x)\, dx = \int_{x_0}^\infty \xi(x)\, dx$ and therefore the relation cannot be inverted if
$\lambda= \tilde{\beta}= \mathcal{R}_{0,k}$. Finally, we get $\mathcal{R}_0= \lim_{k \to \infty} \mathcal{R}_{0,k}=  \tilde{\beta}$ which in this case is independent of the size-dependent growth and diffusion processes.

Notice that \eqref{brnsize} gives the same $\mathcal{R}_0$ under hypothesis \eqref{birth}:
$$
\mathcal{R}_0= \int_{x_0}^\infty \frac{\beta(x)}{\gamma(x)} \, e^{-\int_{x_0}^x \frac{\mu(y)}{\gamma(y)}dy}\, dx=  \tilde{\beta}\; .
$$

\subsection{Age structured model with diffusion in age}\label{age}

In this section we consider a linear age-structured population model for the dynamics of a closed
population of individuals (animals, cells, ...) structured by
\textit{physiological age} $a \geq 0$ with diffusion $D \geq 0$:

\begin{equation}\label{sys}
\left\{ \begin{array}{l}
    \partial_t u(a,t) + \partial_a u(a,t) + \mu\, u(a,t) = D\, \partial_{aa} u(a,t), \smallskip\\
    u(0,t) - D \, \partial_a u(0, t) = \int_0^\infty \beta(a) u(a,t) \, da. 
  \end{array}\right.
\end{equation}
We interpret physiological age as an internal variable associated to
the physiological development of an individual, normally distributed
and with expected value equal to the chronological age.\\
The
physiological (or biological) age of an individual, as opposed to
chronological age or time since birth, varies in a diffusive way. For papers on diffusion age-structured problems see
\cite{Brewer},\cite{Kakumani}. Notice that \eqref{sys} is a particular case of \eqref{linearized}.

The main goal is to compute the basic reproduction number $\mathcal{R}_0^D$
using \eqref{R0} for model (\ref{sys}) and to compare it with the classical
age-structured population model corresponding to $D=0$.

Looking for exponential solutions of system (\ref{sys}) one obtains a characteristic equation from which the following threshold can be derived \cite{Magal}:

\begin{equation} \label{conjecture}
\mathcal{\tilde{R}}_0^D = \frac{2}{1+\sqrt{1+4D \mu}} \int_0^\infty \beta(a)
\, e^{(1-\sqrt{1+4D \mu})\frac{a}{2D} } \, da \; .
\end{equation}
which tends to the
classical expression when the diffusion coefficient tends to zero (by an easy application of the Lebesgue dominated convergence theorem)
$$
\lim_{D \to 0} \mathcal{\tilde{R}}_0^D = \int_0^\infty \displaystyle \beta(a)
\, e^{-\mu a} \, da \; .
$$

This quantity $\mathcal{\tilde{R}}_0^D$ might not be equal to the basic reproduction number
$\mathcal{R}_0^D$ that we will compute using the framework set in Section $2$ but the sign of $\mathcal{\tilde{R}}_0^D-1$ must coincide with the one of $\mathcal{R}_0^D-1$.
In order to obtain $\mathcal{R}_0^D$ we start by computing the inverse of the transition/mortality operator for system (\ref{sys}) which is defined
as $M \phi = \phi' + \mu \phi - D \phi''$. Its inverse can be
computed as the solution $v \in L^1(0,\infty)$ of the boundary value
problem
\begin{equation}\label{sys3}
\left\{ \begin{array}{l}
    v'(a) + \mu\, v(a) - D v''(a) = \phi(a), \smallskip\\
    v(0) = D  v'(0),
  \end{array}\right.
\end{equation}
for a given datum $\phi \in L^1(0, \infty).$ \\
Let us call
$$
\lambda_1 = \frac{1}{2D}(1+\sqrt{1+4D \mu}) > 0 \; , \quad \lambda_2
= \frac{1}{2D}(1-\sqrt{1+4D \mu}) < 0.
$$
The general solution of the second order inhomogeneous linear ode above can be written as:

$$
v(a) = \frac{1}{\sqrt{1+4D \mu}} \left( e^{\lambda_1 a} (c_1 +
\int_a^{\infty} e^{-\lambda_1 s} \phi(s) ds) + e^{\lambda_2 a} (c_2
+ \int_0^{a} e^{-\lambda_2 s} \phi(s) ds) \right)
$$
By imposing that $v$ is integrable we get $c_1=0,$ and using the
boundary condition at $a = 0$ straightforwardly we obtain $c_2=
-\frac{\lambda_2}{\lambda_1} \int_0^{\infty} e^{-\lambda_1 s}
\phi(s) ds$, and thus the solution of (\ref{sys3}) can be written

$$
v(a) = \frac{1}{\sqrt{1+4D \mu}} \left( \int_a^{\infty} e^{\lambda_1
(a-s)} \phi(s) ds + \int_0^{a} e^{\lambda_2 (a-s)} \phi(s) ds
-\frac{\lambda_2}{\lambda_1} \int_0^{\infty} e^{\lambda_2 a
-\lambda_1 s} \phi(s) ds \right) \, .
$$
Therefore we can write

\begin{equation}\label{Mmenys}
(M^{-1} \phi)(a) = \int_0^{\infty} G(a,s) \phi(s) ds,
\end{equation}
where
$$
G(a,s) = \frac{1}{\sqrt{1+4D \mu}} \bigg( e^{\lambda_1(a-s)} H(s-a)
+ e^{\lambda_2(a-s)} H(a-s) -
\frac{\lambda_2}{\lambda_1}e^{\lambda_2 a - \lambda_1 s} \bigg)
$$

Along the lines of \cite{ripoll1} and Section 2 we introduce a
sequence of ``type I" models which approximate (\ref{sys}):

\begin{equation}\label{sysapprox}
\left\{ \begin{array}{l}
    \partial_t u(a,t) + \partial_a u(a,t) + \mu\, u(a,t) = D\, \partial_{aa} u(a,t) +  \int_0^\infty \beta(a) u(a,t) \, da \,k \1_{[0,1/k]} (a) \smallskip\\
    u(0,t) - D \, \partial_a u(0, t) = 0.
  \end{array}\right.
\end{equation}

Notice that \eqref{sysapprox} has a biological meaning on its own.
That is, it corresponds to assuming that the newborns size (or
physiological age) is distributed uniformly on a small interval
around the minimal size $0.$ Similarly to Subsection 3.2, this is an assumption at least as reasonable as assuming that the size of the newborns is concentrated at a given value, for instance $0$.

Here the birth operator $B_k u:=  \int_0^\infty \beta(a) u(a,t) \,
da \,k \1_{[0,1/k]}$ has a one-dimensional range spanned by the
characteristic function of the interval $[0,1/k].$ So, $Lv := \int_0^{\infty} \beta(a) v(a) da$ is a bounded linear form given by the $L^{\infty}$ function $\beta.$ Consequently, the
basic reproduction number for \eqref{sysapprox} can be computed as
$$
\mathcal{R}_{0,k} = \rho \big(B_k M^{-1}\big) = LM^{-1}(k \1_{[0,1/k]}) = \int_0^{\infty} \beta(a) M^{-1}
(k \1_{[0,1/k]})(a) \, da.
$$

In this case, a straightforward computation using \eqref{Mmenys}
gives
\begin{equation}
\label{mort}
\begin{array}{l}
\big(M^{-1} (k \1_{[0,1/k]})\big) (a) = \\
= \frac{1}{\sqrt{1+4 \mu D}} \left[(1-ak)\, \left(F(\lambda_1
(\frac{1}{k}-a)) - F(\lambda_2 (\frac{1}{k}-a)) \right) \,
\1_{[0,1/k]}(a) + \left( F(\frac{\lambda_2}{k}) -
\frac{\lambda_2}{\lambda_1} F(\frac{\lambda_1}{k}) \right) \,
e^{\lambda_2 a} \right]
 \end{array}
\end{equation}
where $F(x):=\frac{1-e^{-x}}{x}$ for $x \neq 0, F(0) = 1.$ \\

We now proceed to compute the basic reproduction number for the limit model \eqref{sys}. So, let us compute the $L^1-$limit of the sequence $M^{-1}(k \1_{[0,1/k]}).$ First notice that $F(x)$ is an analytic function
and that its first derivative is negative and its second derivative is positive.\\
Since $F$ is decreasing and $\lambda_2 < \lambda_1$, the coefficient
of $\1_{[0,1/k]}$ in \eqref{mort} is non-positive, implying
$$\big(M^{-1} (k \1_{[0,1/k]})\big)(a) \leq
\frac{1}{\sqrt{1+4 \mu D}} \left( F(\frac{\lambda_2}{k}) -
\frac{\lambda_2}{\lambda_1} F(\frac{\lambda_1}{k}) \right) \,
e^{\lambda_2 a}.
$$
On the other hand, let us consider the function $G(x) =
F(\frac{\lambda_2}{\lambda_1}x)-\frac{\lambda_2}{\lambda_1} F(x),$
whose derivative $G'(x) = \frac{\lambda_2}{\lambda_1}
\big(F'(\frac{\lambda_2}{\lambda_1} x) - F'(x) \big)$ is positive
for $x>0$ since $\frac{\lambda_2}{\lambda_1} < 0$ and $F''(x) > 0.$
Therefore, $F(\frac{\lambda_2}{k}) - \frac{\lambda_2}{\lambda_1}
F(\frac{\lambda_1}{k}) = G(\frac{\lambda_1}{k})$ decreases with $k$.
So we obtain the following bound, independent of $k$,
$$
\big(M^{-1} (k \1_{[0,1/k]})\big) (a) \leq \frac{1}{\sqrt{1+4 \mu
D}} \big( F(\lambda_2) - \frac{\lambda_2}{\lambda_1} F(\lambda_1)
\big) e^{\lambda_2 a}.
$$
$\big(M^{-1} (k \1_{[0,1/k]})\big) (a)$ converges pointwise to
$\psi_{\infty}(a) := \frac{(1-\frac{\lambda_2}{\lambda_1}) e^{\lambda_2 a}}{\sqrt{1+ 4
\mu D}}$
 and hence the convergence is also in $L^1$ by the Lebesgue dominated
 convergence theorem. Notice that the limit function does not belong
 to the domain of the operator $-M$ because it does not satisfy the
 boundary condition. Also notice that $\psi_{\infty}(a) = G(a,0).$
Along the lines of Section 2, we compute the basic reproduction
number of system \eqref{sys} as
$$\begin{array}{l}
\mathcal{R}_0^D = \lim \mathcal{R}_{0,k} = L\psi_{\infty} = L G(., 0) \\
= \int_0^{\infty} \beta(a) \frac{(1-\frac{\lambda_2}{\lambda_1})
e^{\lambda_2 a}}{\sqrt{1+ 4 \mu D}} \, da = \frac{2}{1+\sqrt{1+4D
\mu}} \int_0^\infty \beta(a) \, e^{(1-\sqrt{1+4D \mu})\frac{a}{2D} }
\, da \; ,
\end{array}
$$
which is the basic reproduction number for the diffusion
age-structured population model (\ref{sys}) and coincides with the
conjecture given at the beginning of the subsection. For an alternative computation of $\mathcal{R}_0^D$ using the solution semigroup we refer the reader to the appendix.
For a constant fertility rate, i.e. $\beta(a) \equiv \bar{\beta}$,
the basic reproduction number is given by
$$
\mathcal{R}_0^D = \frac{\bar{\beta}}{\mu} \; , \quad \textrm{for all $D\geq
0,$}
$$
so it is independent
of the diffusion process.

For age specific fertilities, an interesting issue is the dependence of $\mathcal{R}_0^D$ on the diffusion coefficient $D.$ As expectable, $\mathcal{R}_0^D$ does not depend on $D$ if the fertility function is a constant. On the other hand, we already know that it tends to the basic reproduction number without diffusion when $D$ tends to 0. It is also very easy to show that it tends to 0 for $D$ tending to $\infty$ whenever the fertility function $\beta$ is integrable (an expectable result too), albeit it does it slowly (as $\frac{1}{\sqrt{D}}$).\\
However this does not mean that $\mathcal{R}_0^D$ is always a monotonous (decreasing) function of $D$. For instance the particular case $\beta(a) = \beta_0 a^2 e^{-a}$ (see \cite{Webb} for a discussion on the choice of this birth function) leads to $\mathcal{R}_0^D = \frac{32 \beta_0 D^3}{(2D +\sqrt{1+4 \mu D}-1)^3 (1+\sqrt{1+4 \mu D})},$ which tends to $\frac{2 \beta_0}{(\mu + 1)^3}$ when $D$ goes to $0$ (the value of $\mathcal{R}_0$ for $D = 0$) but, for $\mu>1/2,$ it has a maximum value $\frac{8 \beta_0}{27 \mu}$ (greater than  $\frac{2 \beta_0}{(\mu + 1)^3}$) at $D = 4 \mu - 2.$\\
On the other hand, a very small fertility function in early ages like the one above leads to think in the case of a maturation age below which the individuals do not reproduce. For instance, if $\beta(a)$ vanishes for $a<1$ and it is a constant $\beta_0$ for $a>1$, then $\mathcal{R}_0^D = \frac{\beta_0}{\mu} e^{\frac{1-\sqrt{1+4 \mu D}}{2 D}}$ is an increasing function of $D$ ranging from $\frac{\beta_0}{\mu} e^{- \mu}$ at $D = 0$ to $\frac{\beta_0}{\mu}$ (the basic reproduction number without maturation delay) at infinity.
Then, the presence of a moderate diffusion in age tends to increase $\mathcal{R}_0$ and hence the survival chances of a population.

\begin{figure}[h!]	
\includegraphics[width=0.8\textwidth]{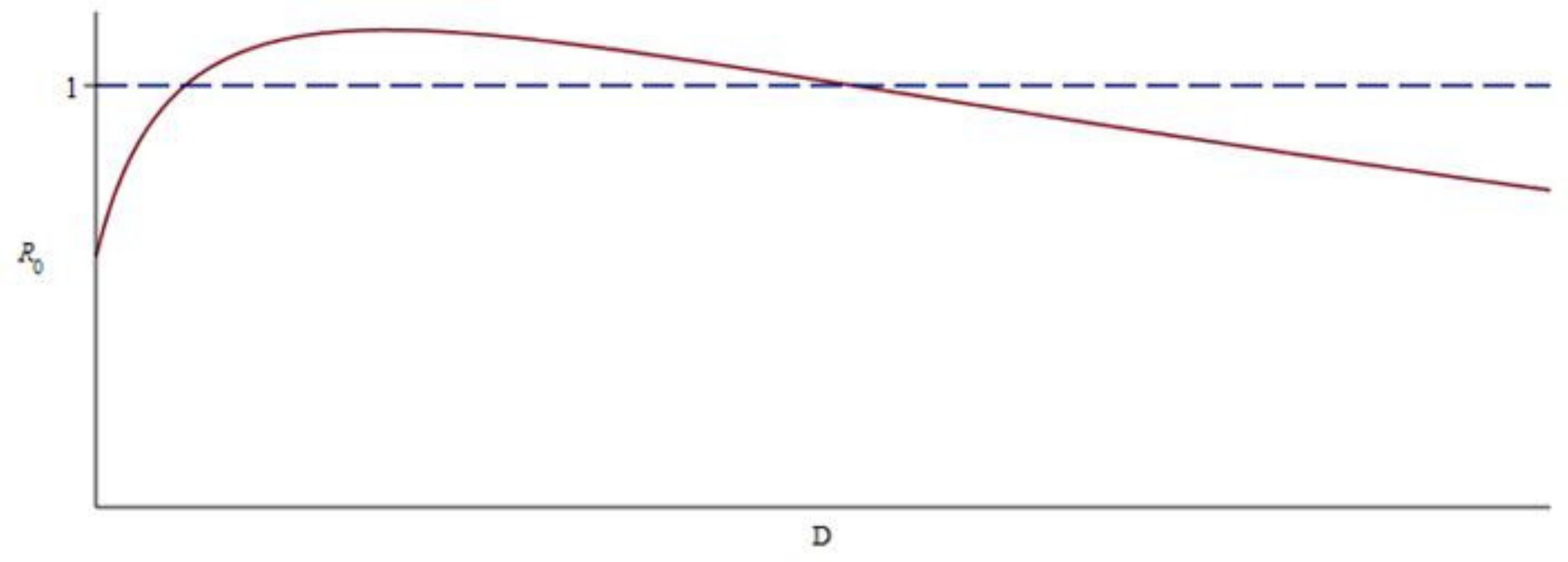}
\caption{A graph of $\mathcal{R}_0^D$ showing values larger than $1$ for moderate diffusion coefficient.}
\end{figure}

\section*{Acknowledgements}

This work has been partially supported by the project
MTM2017-84214-C2-2-P of the Spanish government.

\section{Appendix: An alternative computation of the basic reproduction number for the age structured model with diffusion using the solution semigroup}

An explicit expression for the solution semigroup of the system
\begin{equation}\label{sys2}
\left\{ \begin{array}{l}
    \partial_t v(a,t) + \partial_a v(a,t) + \mu\, v(a,t) = D\, \partial_{aa} v(a,t), \smallskip\\
    v(0,t) - D \, \partial_a v(0, t) = 0, \smallskip\\
    v(a,0) = v_0(a)
  \end{array}\right.
\end{equation}
which is the initial value problem for model \eqref{sys} with zero fertility rate, is possible by means of the Green's function.

First notice that this initial boundary value problem can be rewritten as an initial boundary problem for the heat equation satisfying a slightly different Robin type boundary condition by means of the change of dependent variable $v(a,t) = e^{\frac{a}{2 D} - (\mu + \frac{1}{4D}) t } u(a,t).$ Indeed, a straightforward computation shows that $v$ is a solution of (\ref{sys2}) if and only if $u$ solves

\begin{equation}\label{sys21}
\left\{ \begin{array}{l}
\partial_t u(a,t) = D\, \partial_{aa} u(a,t), \smallskip\\
u(0,t) - 2 D \, \partial_a u(0, t) = 0, \smallskip\\
u(a,0) = e^{-\frac{a}{2D}} v_0(a)
\end{array}\right.
\end{equation}

A Green's function for problem \eqref{sys21} can be found in the book \cite{Cole} (page 602), cataloged as $X30:$

$$
G_{X30}(a,s,t) =
\frac{e^{-\frac{(a-s)^2}{4Dt}}+e^{-\frac{(a+s)^2}{4Dt}}}{2\sqrt{\pi D t}}
-\frac{e^{\frac{t}{4D}+\frac{a+s}{2D}}}{2D}\left(1-\erf
\left(\frac{a+s+t}{2\sqrt{Dt}}\right) \right).
$$
So, the solution of \eqref{sys21} is given by
$$
u(a,t) = \int_0^{\infty} G_{X30} (a,s,t) e^{-\frac{s}{2D}} v_0(s) ds,
$$
and consequently, that of \eqref{sys2} by

$$
(T(t)v_0)(a) := v(a,t) = \int_0^{\infty} G_0(a,s,t)  v_0(s) ds
$$
where

$$\begin{array}{l}
G_0(a,s,t) = e^{\frac{a}{2 D} - (\mu + \frac{1}{4D}) t } e^{-\frac{s}{2D}} G_{X30} (a,s,t)\\
= e^{-\mu t} \left( \frac{e^{\frac{2(a-s)-t}{4D}}}{\sqrt{\pi D t}} \frac{e^{-\frac{(a-s)^2}{4Dt}}+e^{-\frac{(a+s)^2}{4Dt}}}{2} - \frac{e^{\frac{a}{D}}}{2D} \left(1-\erf
\left(\frac{a+s+t}{2\sqrt{Dt}}\right)\right)\right)\end{array}
$$

With the aim of writing the next generation operator, after interchanging the integration order we can write
$$
\int_0^{\infty} (T(t) v)(a) dt = \int_0^{\infty} \int_0^{\infty}  G_0(a,s,t)\, dt \, v(s) \, ds
$$

Now, proceeding in a similar manner as we did in Section \ref{age}, we consider system \eqref{sysapprox} and, setting $v_k = k \1_{[0,1/k]},$ we can write

$$
\int_0^{\infty} (T(t) v_k)(a) dt = k \int_0^{1/k} \int_0^{\infty}  G_0(a,s,t)\, dt \, ds \rightarrow  \int_0^{\infty}  G_0(a,0,t)\, dt,
$$

where
$$
G_0(a,0,t) = e^{-\mu t} \cdot \left(
\frac{e^{-\frac{(a-t)^2}{4Dt}}}{\sqrt{\pi D t}}
-\frac{e^{\frac{a}{D}}}{2D}\left(1-\erf
\left(\frac{a+t}{2\sqrt{Dt}}\right) \right) \right).
$$

Arguing as above,

$$
\mathcal{R}_0^D = \lim \mathcal{R}_{0,k} = \lim \int_0^{\infty} \beta(a) \, (T(t) v_k)(a) \, dt = \int_0^{\infty} \beta(a) \int_0^{\infty}  G_0(a,0,t)\, dt \, da,
$$

and the only thing left is the computation of the integral of $G_0(a,0,t)$ with respect to $t$ on $(0, \infty)$.

Indeed, a standard argument based on ``completing squares" gives

\begin{equation}\label{integral}
\int_0^\infty \, e^{-\mu \tau} \cdot
\frac{e^{-\frac{(a-\tau)^2}{4D\tau}}}{\sqrt{\pi D \tau}} d\tau =
2\frac{e^{\frac{a}{2D}(1-\sqrt{1+4D\mu})}}{\sqrt{1+4D \mu}},
\end{equation}
whereas for the second term we have,

$$
\int_0^\infty e^{- \mu \tau} \frac{e^{\frac{a}{D}}}{2D}\left(1- \erf
\left(\frac{a+\tau}{2\sqrt{D\tau}}\right) \right) d\tau \; =
\frac{e^{\frac{a}{D}}}{2D}\int_0^\infty e^{-\mu \tau}
\int_{\frac{a+\tau}{2\sqrt{D \tau}}}^\infty
\frac{2}{\sqrt{\pi}}e^{-x^2} \, dx \, d \tau
$$
which, after changing integration order equals

$$ \frac{e^{\frac{a}{D}}}{D \sqrt{\pi}} \int_{\sqrt{a/D}}^\infty
e^{-x^2}\int_{2Dx^2-a-2Dx\sqrt{x^2-\frac{a}{D}}}^{2Dx^2-a+2Dx\sqrt{x^2-\frac{a}{D}}}
e^{-\mu \tau} \, d\tau \, dx $$

$$ = \frac{e^{(1+\mu D)\frac{a}{D}}}{\mu D \sqrt{\pi}}\left(
\int_{\sqrt{a/D}}^\infty e^{-(1+2 \mu D)x^2 + 2 \mu D x
\sqrt{x^2-\frac{a}{D}}} \, dx - \int_{\sqrt{a/D}}^\infty e^{-(1+2
\mu D)x^2 - 2 \mu D x \sqrt{x^2-\frac{a}{D}}} \, dx \right)$$
$$
= \frac{e^{(1+\mu D)\frac{a}{D}}}{\mu D \sqrt{\pi}}
\left(\int_{\sqrt{\frac{a}{D}}}^{\infty} + \int_0^{\sqrt{\frac{a}{D}}} \right) e^{-\frac{(s^2+
\frac{a}{D})(s^2 + (4\mu D+1)\frac{a}{D})}{4s^2}} \big( \frac{1}{2}
- \frac{a/D}{2s^2} \big) \, ds,
$$
where we replaced $x=\frac{1}{2} (s+\frac{a/D}{s}), s > \sqrt{\frac{a}{d}}$ in the first integral
and $x=\frac{1}{2} (s+\frac{a/D}{s}), s\in (0,\sqrt{\frac{a}{d}})$ in the second one. This finally
gives
$$
\frac{e^{(1+\mu D)\frac{a}{D}}}{2 \mu D \sqrt{\pi}} \left( \int_0^{\infty}  e^{-\frac{(s^2+
\frac{a}{D})(s^2 + (4\mu D+1)\frac{a}{D})}{4s^2}}  \, ds -  \int_0^{\infty} \frac{a}{D s^2} e^{-\frac{(s^2+
\frac{a}{D})(s^2 + (4\mu D+1)\frac{a}{D})}{4s^2}}  \, ds \right)
$$
$$\begin{array}{l}
= \frac{e^{\frac{a}{2 D}}}{4 \mu D}  \biggl[ e^{\frac{a \sqrt{1 + 4 \mu D}}{2 D}} \erf \big( \frac{s}{2} +\frac{a \sqrt{1 + 4 \mu D}}{2 D s} \big) +  e^{-\frac{a \sqrt{1 + 4 \mu D}}{2 D}} \erf \big( \frac{s}{2} -\frac{a \sqrt{1 + 4 \mu D}}{2 D s}  \big) \Big|_0^{\infty}  \\
 + \frac{1}{\sqrt{1 + 4 \mu D}} \left( e^{\frac{a \sqrt{1 + 4 \mu D}}{2 D}} \erf \big( \frac{s}{2} +\frac{a \sqrt{1 + 4 \mu D}}{2 D s} \big) -  e^{-\frac{a \sqrt{1 + 4 \mu D}}{2 D}} \erf \big( \frac{s}{2} -\frac{a \sqrt{1 + 4 \mu D}}{2 D s} \big) \right) \Big|_0^{\infty} \bigg]
\end{array}
$$

$$
= \frac{e^{\frac{a}{2 D}(1- \sqrt{1 + 4 \mu D})}}{2 \mu D} \big( 1 -
\frac{1}{\sqrt{1 + 4 \mu D}} \big) ,
$$
which can be checked by direct differentiation. \\
Subtracting the last expression from (\ref{integral}) one easily
obtains (\ref{conjecture}).

%
%

\end{document}